\documentclass{amsart}
\usepackage{amssymb}
\newcommand{\bor}{{\mathsf   {Borel}}}
\newcommand{\Cohen}{{\mathbb C}}
\newcommand{\xp}{{\mathbf X}}
\newcommand{\PM}{{\mathbf {PM}}}
\newcommand{\UM}{{\mathbf {UM}}}
\newcommand{\UN}{{\mathbf {UN}}}
\newcommand{\cl}{\operatorname{{\mathsf cl}}}
\newcommand{\intr}{\operatorname{{\mathsf interior}}}

\newcommand{\ZFCa}{{\operatorname{\mathsf {ZFC}}}}

\newcommand{\rest}{{\mathord{\restriction}}}

\newcommand{\N}{{\mathcal N}}
\newcommand{\M}{{\mathcal M}}

\newcommand{\<}{\langle}
\renewcommand{\>}{\rangle}

\newcommand{\Proof}{{\sc Proof} \hspace{0.2in}}
\newcommand{\lft}[2]{\mathopen\ifcase#1{}\oo\or
                        \big#2\or\Big#2\else\oo\fi} 
\newcommand{\rgt}[2]{\mathclose\ifcase#1{}\oo\or
                        \big#2\or\Big#2\else\oo\fi}

\newcommand{\T}{{\mathbf   {T}}}

\theoremstyle{plain}
\newtheorem{theorem}{Theorem}
\theoremstyle{plain}
\newtheorem{lemma}[theorem]{Lemma}

\newtheorem{definition}[theorem]{Definition}

\begin{document}
\title{Remarks on small sets of reals}
\author{Tomek Bartoszynski}
\address{Department of Mathematics and Computer Science\\
Boise State University\\
Boise, Idaho 83725 U.S.A.}
\thanks{Author  partially supported by 
NSF grant DMS 9971282 and Alexander von Humboldt Foundation} 
\email{tomek@math.boisestate.edu, http://math.boisestate.edu/\char 126 tomek}
\begin{abstract}
We show that the Dual Borel Conjecture implies that
${\mathfrak d}> \boldsymbol\aleph_1 $ and find some topological
characterizations of perfectly meager and universally meager sets.
\end{abstract}
\maketitle

\section{Introduction}
For $f,g \in \omega^\omega $ let $f \leq^\star g$ mean that $f(n) \leq 
g(n)$ for all but finitely many $n$.
Let 
$${\mathfrak b}=\min\{|F|: F \subseteq \omega^\omega \ \&\ \forall h
\in \omega^\omega \ \exists f \in F \ f \not \leq^\star h\},$$

$${\mathfrak d}=\min\{|F|: F \subseteq \omega^\omega \ \&\ \forall h
\in \omega^\omega \ \exists f \in F \ h  \leq^\star f\}.$$

Let $\N$ be the ideal of measure zero subsets of $2^\omega $ with
respect to the standard product measure $\mu$, and let $\M$ be the
ideal of meager subsets of $2^\omega $.
Let $+$ be the addition mod $2$ on $2^\omega $. For $A,B \subseteq
2^\omega $ let $A+B=\{a+b: a\in A, b\in B\}$.

\begin{definition}
Let $X \subseteq 2^\omega $. We say that:
\begin{enumerate}
\item $X$ has strong measure zero if $X+F\neq 2^\omega $ for all $F
  \in \M$,
\item $X$ is strongly meager if $X+F\neq 2^\omega $ for all $F
  \in \N$, 
\item $X$ is meager additive if
$X+F\in \M $ for all $F
  \in \M$,
\item $X$ is null additive if $X+F\in \N $ for all $F
  \in \N$.
\end{enumerate}
\end{definition}

Let Borel Conjecture be the statement that there are no uncountable
strong measure zero sets, and the Dual Borel Conjecture that there are
no uncountable strongly meager sets.

Rothberger showed that if ${\mathfrak b}=\boldsymbol\aleph_1 $ then
there is an uncountable strong measure zero set. 
We will strengthen this result by showing that:
\begin{theorem}
  \begin{enumerate}
  \item If ${\mathfrak b}=\boldsymbol\aleph_1 $ then there exists an
    uncountable meager additive set.
  \item If ${\mathfrak d}=\boldsymbol\aleph_1 $ then there exists an
    uncountable null additive set. In particular, the Dual Borel
    Conjecture implies that ${\mathfrak d}> \boldsymbol\aleph_1 $
  \end{enumerate}
\end{theorem}
\begin{proof}
The construction presented here is a modification of a construction
invented by Todorcevic.
Part (1) was also proved in \cite{GalMil84Gam} using different methods.
  
Let $F=\{f_\alpha:\alpha<\omega_1\}$ be a family of functions in
$\omega^\omega $ such that 
\begin{enumerate}
\item $f_\alpha $ is strictly increasing for $\alpha<\omega_1$,
\item $\forall \alpha<\beta \ f_\alpha \leq^\star f_\beta $.
\end{enumerate}

For a perfect tree $p \subseteq 2^{<\omega}$ let $[p]$ denote the set
of its branches. Every perfect subset of $2^\omega $ is a set of
branches of a perfect tree.

We will build an $\omega_1$-tree $\T$ of perfect subsets of $2^\omega
$. Let $\T_\alpha $ denote the $\alpha $-th level of $\T$.
We require that
\begin{enumerate}
\item $\forall \beta>\alpha \ \forall n \in \omega \ \forall p \in
  \T_\alpha \ \exists q\in T_\beta \ \left(q \subseteq p \ \&\ q \cap
    2^n=p\cap 2^n\right)$,
\item $\forall p \in \T_{\alpha+1} \ \forall^\infty n \ |p \cap
  2^{f_\alpha(n)}| \leq 2^n$,
\item $\forall \alpha \ \T_\alpha $ is countable.
\end{enumerate}

{\sc Successor step}.
Suppose that $\T_\alpha $ is given. For each $p \in \T_\alpha $ choose
$\{q_n: n \in \omega\}$ such that 
\begin{enumerate}
\item $ \forall n \ q_n \subseteq p$,
\item $q_n \cap 2^n=p\cap 2^n$,
\item $[q_n] \cap [q_m]=\emptyset$ for $n\neq n$,
\item $\forall n\ \forall^\infty k\ |q_n \cap
  2^{f_\alpha(k)}| \leq 2^k$.
\end{enumerate}

Set $\{q_n: n \in \omega \}$ to be the successors of $p$ on level
$\T_{\alpha+1}$.

{\sc Limit step} Suppose that $\gamma $ is  a limit ordinal and
$\{\T_\alpha: \alpha < \gamma\}$ are already constructed.

For each $p \in \bigcup_{\alpha<\gamma}\T_\alpha $ and $n \in \omega $
we will construct an element $q=q(p,n)$ belonging to the level
$\T_\gamma $ as follows.

Suppose $p=p_0 \in \T_{\alpha_0}$, $n=n_0$ and construct sequences $\<\alpha_k:
k\in \omega \>$, $\<n_k: k \in \omega \>$ and $\<p_k: k\in \omega \>$
such that 
\begin{enumerate}
\item $p_k \in \T_{\alpha_k}$,
\item $\sup_k \alpha_k=\gamma $, $\lim_k n_k =\infty$,
\item $p_{k+1} \subseteq p_k$ for all $k$,
\item $p_{k+1} \cap 2^{n_k}=p_k \cap 2^{n_k}$,
\item $q = \bigcap_k p_k$ is a perfect tree.
\end{enumerate}
The last condition is guaranteed by the careful choice of the sequence
$\<n_k: k\in \omega \>$.

Let $X$ be the set obtained by selecting one element out of every tree
$p \in \T$.

The following lemma gives the first part of the theorem.
\begin{lemma}
  If $F$ is an unbounded family in $\omega^\omega $ then $X$ is meager
  additive.
\end{lemma}
\Proof
Suppose that $H \subseteq 2^\omega $ is a meager set. It is well known
(see \cite{BJbook}, theorem 2.2.4) that 
there exists $x_H \in
  2^\omega$ and a strictly increasing function $f_H \in \omega^\omega$
  such that 
$$H \subseteq  \{x \in 2^\omega : \forall^\infty n \ \exists j \in
\lft1[f_H(n), f_H(n+1)\rgt1) \ x(j) \neq x_H(j)\}.$$
Since we work with translations without loss of generality we can
assume that $x_H(k)=0$ for all $k$.
As $F$ is unbounded there is $\alpha_0<\omega_1$ such that 
$$\exists^\infty n \ \exists k\ f_{\alpha_0}(n)<f_H(k)<f_H(k+1)<\cdots
< f_H(k+2^n)<f_{\alpha_0}(n+1).$$
Fix sequences $\<u_n,k_n: n \in \omega \>$ such that 
$$\forall  n \ \ f_{\alpha_0}(u_n)<f_H(k_n)<f_H(k_n+1)<\cdots
< f_H(k_n+2^{u_n})<f_{\alpha_0}(u_n+1).$$
Fix $p \in \T_{\alpha_0+1}$ and let $z_p\in 2^\omega $ be defined as
follows:
given $n \in \omega $ let $\{s_1, \dots, s_{2^{u_n}}\}$ be an
enumeration of $p \cap 2^{f_{\alpha_0}(u_n+1)}$.
Define
$z_p \rest [f_H(k_n+i),f_H(k_n+i+1))=s_i\rest
[f_H(k_n+i),f_H(k_n+i+1))$ for $i \leq 2^{u_n}$.
This defines $z_p$ on an infinite subset of $\omega $, extend it
arbitrarily to a total function.
Let
$$G=\{x \in 2^\omega : \forall^\infty n \ x \rest
[f_{\alpha_0}(u_n),f_{\alpha_0}(u_n+1))\neq z_p
[f_{\alpha_0}(u_n),f_{\alpha_0}(u_n+1))\}.$$
We claim that $[p]+H \subseteq G$.
Note that if $x \in [p]+H$ then there exists $y \in [p]$ such that 
$$\forall^\infty n \ x \rest [f_H(n),f_H(n+1))\neq y\rest
[f_H(n),f_H(n+1)).$$
Fix one such $y$ and note that
for sufficiently large $n$ there is $i$ such that
$y \rest [f_H(k_n+i),f_H(k_n+i+1))=s_i\rest
[f_H(k_n+i),f_H(k_n+i+1))=z_p \rest [f_H(k_n+i),f_H(k_n+i+1)),$
which implies that
$x \rest [f_H(k_n+i),f_H(k_n+i+1))\neq z_p \rest
[f_H(k_n+i),f_H(k_n+i+1))$.
It follows that
$x\rest [f_{\alpha_0}(u_n),f_{\alpha_0}(u_n+1)) \neq z_p \rest
[f_{\alpha_0}(u_n),f_{\alpha_0}(u_n+1)),$
which means that $x \in G$.

Let $X_{\alpha_0}$ be the collection of points selected from levels
$\bigcup_{\alpha<\alpha_0}\T_\alpha $. We have
$$X+G \subseteq (X_{\alpha_0}+G)  \cup \bigcup_{p \in \T_{\alpha_0}}
[p]+G \in \M,$$ which finishes the proof.

To prove the second part of the theorem we will use the following
lemma:
\begin{lemma}
  If $F$ is a dominating family then $X$ is null additive.
\end{lemma}
\Proof
The following is well known:
\begin{lemma}
  Suppose that $H \subseteq 2^\omega $ is a null set. 
There
  exists a sequence of clopen sets $\{C_n:n \in \omega\}$ such that 
  that for all $n$,
  \begin{enumerate}
  \item $\mu(C_n)< 4^{-n}$,
  \item $ \forall x \in H \ \exists^\infty n \ x  \in C_n$.
  \end{enumerate}
\end{lemma}
\begin{proof}
Since $H$ has measure zero there are open sets
$\<U_{n} : n \in \omega \>$ 
covering $H$ such that $\mu(U_{n}) <  4^{-n-1}$,
for $n \in \omega$.
Write each set $U_{n}$ as a disjoint union of open basic intervals,
$U_{n} = \bigcup_{m=1}^{\infty} [s^{n}_{m}]$
for  $n \in \omega$, and 
order these sequences lexicographically in a single sequence $\{t_n: n
\in \omega\}$.
Put $k_0=0$ and for $n>0$ let 
$k_{n}=\min\{k: \sum_{k>k_{n}} \mu([t_k])<
4^{-n-1}\}$.
Let 
$$C_{n}=\bigcup_{k \in [k_n,k_{n+1})} [t_k].$$
 
Clearly $\mu(C_n)< 4^{-n}$ and, since all basic sets have been accounted for,
$$ \forall x \in H \ \exists^\infty n \ x  \in C_n.$$
\end{proof}

Let $H \subseteq 2^\omega $ be a measure zero set and $\{C_n: n \in
\omega\}$ the sequence given by the lemma.
Since each clopen set is a union of finitely many basic sets we can
find a function $f_H \in \omega^\omega $ such that for every $n$, $C_n
\subseteq 2^{f_H(n)}$.
Let $\alpha_0$ and $n_0$ be such that 
$$\forall n>n_0 \ f_H(n)<f_{\alpha_0}(n).$$
By modifying finitely many $C_n$'s we can assume that $n_0=0$.
Suppose that $p \in \T_{\alpha_0+1}$. 
Note that
$$[p]+C_n \subseteq (p \cap 2^{f_H(n)})+C_n=D_n.$$
Since $|p \cap 2^{f_H(n)}|\leq |p \cap 2^{f_{\alpha_0}(n)}|$, it
follows that
$$\mu\left([p]+C_n\right)=\mu(D_n) \leq 2^n \cdot 4^{-n} \leq
2^{-n}.$$
Therefore, $[p]+H \subseteq G=\{x \in 2^\omega : \exists^\infty n \
x \in D_n\}$, and $\mu(G)=0$. 
The rest of the proof is identical to the proof of the first part.
\end{proof}

\section{Perfectly meager and universally meager sets}
A set $X \subseteq 2^\omega $ is perfectly meager ($X \in \PM$) if $P \cap X$ is
meager in $P$ for every perfect set $P$.

A set $X \subseteq 2^\omega $ is universally meager ($X \in \UM$,
\cite{Zakr}) if for every Borel isomorphism $F : 2^\omega
\longrightarrow 2^\omega $, $F"(X) \in \M$.

A set $X \subseteq 2^\omega $ is universally null ($X \in \UN$) if for
every Borel isomorphism $F : 2^\omega 
\longrightarrow 2^\omega $, $F"(X) \in \N$.

In \cite{Zakr} it is shown that the
universally meager sets are a category analog of universally null sets.
Yet, for quite a while perfectly meager sets were viewed in this
role. Theorems \ref{lutzer} and \ref{two} explain this phenomenon, which has to
do with the fact that the families $\PM$ and $\UM$ are not {\em very} different.
It is easy to see that $\UM \subseteq \PM$. Continuum hypothesis
 or Martin's Axiom imply that
$\UM\neq \PM$, but it is also consistent that $\UM=\PM$, \cite{BaPerf00}.

The first theorem gives a simple proof of a characterization of
perfectly meager sets found in \cite{bhm99}. 
 
\begin{theorem}\label{lutzer}
The following are equivalent:
\begin{enumerate}
\item $X \in \PM$,
\item for every countable dense-in-itself set $A$ there exists a set $B \subseteq A$
  such that $\cl(A)=\cl(B)$ and $B$ is a $G_\delta $-set relative to
  $X\cup A$,
\item  for every countable  set $A$ there exists a set $B \subseteq A$
  such that $\cl(A)=\cl(B)$ and $B$ is a $G_\delta $-set relative to
  $X\cup A$,
\item there exists
an $F_\sigma $ set $F$ such that $X \subseteq F$ and $F$ meager in $P$.
\end{enumerate}
\end{theorem}
\begin{proof}
(2) $ \rightarrow $ (1).
Let $P$ be a perfect set and $Q \subseteq P$ a countable dense set. By
(2) without loss of generality we can assume that $Q$ is a $G_\delta $
relative to $X$. In other words $Q=X \cap \bigcap G_n$, where $G_n$
are open sets. Clearly $G_n$'s are dense in $P$. It follows that $X
\cap P \subseteq Q \cup P \setminus \bigcap G_n$.

\bigskip 
(1) $\rightarrow $ (4)
Since $X \in \PM$ there exists an $F_\sigma $ set $F_1$ such that $X
\subseteq F_1$, and there exists an $F_\sigma $ set $F_2$ such that
$F_2 \cap P$ is meager in $P$ and $X \cap P \subseteq F_2$.
Now the set $F = (F_1 \setminus P) \cup (F_2 \cap P)$ is
the set we are looking for.

\bigskip

(4) $\rightarrow $ (3).
Suppose that $A$ is countable. By Cantor-Bendixson Theorem,
\cite{Kechris95}, $\cl(A)=P \dot{\cup} C$, where $P$ is perfect, $C$
is countable and open relative to $P$.
If $P=\emptyset$ then $C$ is countable and closed and $A=C \setminus
(C\setminus A)$ is a $G_\delta $ set in $2^\omega $.
Thus assume that $P \neq\emptyset$
and let $\<F'_n: n \in \omega\>$ be closed nowhere dense sets  such
that $X  \subseteq \bigcup_n F'_n$ and $F'_n \cap P$ is closed nowhere
dense in $P$ for each $n$. 
Let $\{F_n: n \in \omega\}$ be closed sets such that 
$\bigcup_n F_n = \left(\bigcup_{n \in \omega} F'_n\right) \setminus
C$. Let $A'=A \setminus C$ and consider sets 
$A_n = A' \setminus F_n$. Since the family $\<A_n: n \in \omega\>$
has the finite intersection property we can find a set $B' \subseteq A'$
such that 
\begin{enumerate}
\item $B' \setminus A_n $ is finite for all $n$,
\item $B'$ is dense in $P$.
\end{enumerate}
Let 
$F_n^\star = F_n \setminus B'$. Note that each set $F_n^\star $ is
$F_\sigma $
since it differs from $F_n$ by a finite set. 
Put
$B = (X\cup A) \cap \bigcap_n (2^\omega  \setminus F^\star_n).$ It
follows that $B' \cup (A \cap C) \subseteq B \subseteq A$, which finishes the proof.

\bigskip 

(3) $\rightarrow $ (2) is obvious.
\end{proof}

\begin{theorem}\label{two}
The following are equivalent:
\begin{enumerate}
\item $X \in \UM$
\item for every sequence of countable dense-in-itself sets $\{A_n:n
  \in \omega\}$ there exists a sequence  $B_n \subseteq A_n$
  such that $\cl(A_n)=\cl(B_n)$ and $\bigcup_n B_n$ is a $G_\delta
  $-set relative to $X \cup \bigcup_n A_n$.
\item For every sequence of perfect sets $\{P_n: n \in \omega\}$ there
  exists an $F_\sigma $ set $F$ such that $X \subseteq F$ and $F$ is
  meager in $P_n$ for every $n$.
\end{enumerate}
\end{theorem}
\begin{proof}
  (2) $\rightarrow $ (1)
The following argument is a small modification of a proof from \cite{Nowweiss00}.
Suppose that $X \not \in \UM$. Let $G:X \longrightarrow Y'$ be a Borel
isomorphism onto a non-meager set. By Kuratowski's theorem, $G^{-1}$
is continuous on a dense $G_\delta $ set. It follows that there exists a continuous
one-to-one
function $F:Y \longrightarrow X$, where $Y \not\in \M$.
Let $\{U_n: n \in \omega\}$ be enumeration of clopen subsets of
$2^\omega $ such that $U_n \cap Y$ is uncountable.
For each $n$, choose a countable dense-in-itself set $A_n \subseteq
F"(Y \cap U_n)$.
We will show that every $G_\delta $ set  which is disjoint from $X \setminus
\bigcup_n A_n$ is also disjoint from one of the $A_n$'s.
Let $F=\bigcup_n F_n$ be an $F_\sigma $ set containing $X \setminus
\bigcup_n A_n$. For every $n$ let $H_n$ be a closed set such that
$H_n \cap Y=F^{-1}(F_n)$.
If for every $n$, $Y_n=\intr(H_n) \cap Y$ is countable then 
$Y \subseteq \bigcup_n Y_n \cup \bigcup_n \left(H_n\setminus \intr(H_n)\right)\in
\M$, which is impossible. Thus there exists $m, n \in \omega $ such that $A_m
\subseteq F"(U_m \cap Y) \subseteq  F_n$.

\bigskip

(1) $\rightarrow $ (3)

Let $\Cohen$ denote the Cohen algebra. 
The following is a (small) fragment of the theorem 2.1 from \cite{Zakr}.
\begin{theorem}\label{zak} 
For a subset $X$ of a perfect Polish space $\xp$, the
following are equivalent:
\begin{enumerate}
\item $X \in \UM$.
\item For every $\sigma$-ideal $ {\mathcal J}$ in $\bor(\xp)$ such that
  $\bor(\xp)/ {\mathcal J} 
\cong \Cohen$ there is a Borel set $B\in {\mathcal J} $ such that $X\subseteq
B$.
\end{enumerate}
 \end{theorem} 
 \begin{proof}
   The implication (1) $\rightarrow $ (2) that is of interest to us is
   a consequence of Sikorski's theorem (\cite{Kechris95}, 15.10): if
   ${\mathcal J} $ is a
$\sigma $-ideal in $\bor(\xp)$ such that $\bor(\xp)/ {\mathcal J}  \cong \Cohen$,
then there is Borel automorphism $F: \xp \longrightarrow \xp$ such
that 
$$
\forall X\in\bor(\xp) \ X\in\M \iff F"(X)\in {\mathcal J} . $$
 \end{proof}
Suppose that $\{P_n: n \in \omega\}$ are given. Consider the ideal
$${\mathcal J} = \{A: A \text{ is Borel and } \forall n \ A \cap P_n
\text{ is meager in } P_n\}.$$
It is easy to see that $\bor(\xp)/ {\mathcal J} 
\cong \Cohen$ (as $\bor(\xp)/ {\mathcal J}$ is atomless and has a
countable dense subset). Therefore there exists a set $F \in {\mathcal J} $ such
that $X \subseteq F$.

\bigskip

(3) $\rightarrow $ (2)
Let $P_n=\cl(A_n)$.
Let $F=\bigcup_n F_n$ be an $F_\sigma $ set containing $X$ such that
$F$ is meager in each $P_n$. 
 As in the proof of \ref{lutzer}, 
build by induction a set $B' \subseteq \bigcup_n A_n$ such that for $n
\in \omega $, 
\begin{enumerate}
\item $\cl(B' \cap A_n)=\cl(A_n)$,
\item $B' \setminus F_n$ is finite.
\end{enumerate}
As before $F^\star_n = F_n \setminus B'$ is an $F_\sigma $ set and  let
$B=(X \cup \bigcup_n A_n) \cap \bigcap_n (2^\omega \setminus
F^\star_n)$. The sets $B_n = B \cap A_n$ for $n \in \omega $ are as required.
\end{proof}

The assumption in theorem \ref{two}(2) that the sets $A_n$ are dense-in-itself is
necessary since
we have the following:
\begin{theorem}
  There exists a set $X \in \UM$, and a family $\{A_n:n \in \omega\}$
  of countable subsets of $X$ such that if $B_n \subseteq A_n $ is
  such that $\cl(B_n)=\cl(A_n)$ for each $n$, then $\bigcup_n B_n$ is
  not a $G_\delta $ relative to $X$.
\end{theorem}
\begin{proof}
Suppose that $X \in \UM$, $C \subseteq X$ is countable and $C$ is not
a $G_\delta $ relative to $X$.
Write $C=\bigcup_n A_n$, where each $A_n $ is discrete, that is
$x \not\in \cl(A_n \setminus \{x\})$ for $x \in A_n$. It follows
that if $B_n \subseteq A_n$ and $\cl(B_n)=\cl(A_n)$ then $A_n=B_n$,
which finishes the proof.

To construct a set $X$ as above,
recall that  a set of reals is called
 a $\lambda $-set if all of its countable subsets are relative
 $G_\delta $ sets.
Theorem \ref{two} implies readily that
all $\lambda $-sets, and unions of countable sets with $\lambda$-sets,  are universally meager.

Let $Y \subseteq 2^\omega $ and a countable set $C \subseteq 2^\omega
$ be such that 
\begin{enumerate}
\item $Y$ is a $\lambda $-set,
\item $X=Y \cup C$ is not a $ \lambda $-set, that is $C$ is not a
  $G_\delta $ set relative to $Y$.
\end{enumerate}
Rothberger showed that such sets can be constructed in $\ZFCa$,
(5.6 of \cite{Mil84Spe}).

\end{proof}


\end{document}